\documentclass[12pt,a4paper,draft]{article}
\usepackage[english]{babel}
\usepackage{latexsym,amsfonts,amssymb,amsmath,longtable,amsthm,mathrsfs}
\renewcommand{\leq}{\leqslant}
\renewcommand{\geq}{\geqslant}

\renewcommand{\qed}{\hfill{$\Box$}}

\newtheorem{Lemma}{{\bfseries Lemma}}

\newtheorem{Theo}[Lemma]{{\bfseries Theorem}}
\theoremstyle{definition}
\newtheorem{Conj}{{\bfseries Conjecture}}

\DeclareMathOperator{\PSL}{PSL}
\DeclareMathOperator{\Sym}{Sym} 
 
\DeclareMathOperator{\Alt}{Alt}
 \DeclareMathOperator{\Hall}{Hall}

\newcommand{\splitext}{\,\colon\!}
\newcommand{\arbitraryext}{\,\ldotp}
\newcommand{\nonsplitext}{\,{}^{\text{\normalsize{\textperiodcentered}}}}

\begin{document}

MSC2010 20D20, 20D05, 20D06

\begin{center}

{\bfseries Pronormality of Hall subgroups in finite simple groups\footnote{The work is supported by  RFBR, projects  10-01-00391, 11-01-00456, and
11-01-91158, Federal Target Grant ``Scientific and
educational personnel of innovation Russia'' for 2009-2013 (government
contract  No. 14.740.11.0346).}}

Evgeny P. Vdovin, Danila O. Revin

\end{center}


\begin{center}
Abstract
\end{center}

We prove that Hall subgroups of finite simple groups are pronormal. Thus we obtain an affirmative answer to Problem 17.45(a) of ``Kourovka notebook''.

\section*{Introduction}

According to definition of P.Hall, a subgroup $H$ of a group $G$ is called
{\it pronormal}, if for every $g\in G$ subgroups $H$ and $H^g$ are conjugate in  $\langle H, H^g\rangle$. Classical examples of pronormal subgroups are:

$\bullet$~normal subgroups;

$\bullet$~maximal subgroups;

$\bullet$~Sylow subgroups of finite groups;

$\bullet$~Carter subgroups (i.~e., nilpotent selfnormalizing subgroups) of finite solvable groups;

$\bullet$~Hall subgroups (i.~e. subgroups whose order and index are coprime) of finite solvable groups.

The pronormality of subgroups in last three cases follows from conjugacy of
Sylow, Carter, and Hall subgroups in finite groups in corresponding classes. In
[1, Theorem~9.2] the first author proved that Carter subgroups in finite groups
are conjugate. As a corollary it follows that Carter subgroups of finite groups
are pronormal.

In contrast with Carter subgroups, Hall subgroups in finite groups can
be non-con\-j\-u\-ga\-te. The goal of the authors is to find classes of finite
groups with pronormal Hall subgroups. In the present paper the following result
is obtained.

\begin{Theo}
Hall subgroups of finite simple groups are pronormal.
\end{Theo}

The theorem gives an affirmative answer to Problem 17.45(a)  from the
``Kourovka notebook'' [2],  and it is announced by the authors in [3,
Theorem~7.9]. This result is supposed to use for studying the problem, whether
$C_\pi$ is inherited by overgroups of $\pi$-Hall subgroups  [2,
Problem~17.44(a);
4, Conjecture 3; 5, Problems 2, 3] (all definitions are given below).

\section{Notation, conventions, and preliminary results}

Notation in the paper are standard.

If $G$ is a finite group, $H$ is its subgroup, and $x$ is an element of  $G$,
then by $Z(G)$, $O_\infty(G)$, $N_G(H)$,
$C_G(H)$, and $C_G(x)$ we denote the center of~$G$, the solvable radical
of~$G$, the normalizer of $H$ in $G$, the centralizer of  $H$ in $G$, and the
centralizer of  $x$ in $G$ respectively. Given  groups $A$ and $B$ by
$A\times B$ and  $A\circ B$ we denote the direct product  and a central product
respectively. If $A$ and $B$ are subgroups of $G$, then by  $\langle
A,B\rangle$ and $[A,B]$ the subgroup generated by  ${A\cup B}$ and mutual
commutant of $A$ and $B$ are denoted.

We often use the notations from [6].   In particular, by   $A\,:\, B$,
$A\,\dot{}\,B$, and
$A\, .\, B$ we denote a split, a nonsplit, and an arbitrary extensions of $A$ by
$B$ respectively. Given group $G$ and a subgroup $S$ of the symmetric group
$\operatorname{Sym}_n$ we denote permutation wreath product of $G$ and $S$ by
$G\wr S$ (here $n$ and the embedding of $S$ into $\operatorname{Sym}_n$ assumed
to be known).

We write $H\operatorname{prn} G$ if $H$ is a pronormal subgroup of~$G$.

Throughout  $\pi$ denotes a set of primes. A natural number $n$ with
$\pi(n)\subseteq\pi$, is called a
$\pi$-{\it number}, while a group  $G$ with $\pi(G)\subseteq
\pi$ is called a $\pi$-{\it group}. Symbol $n_\pi $ is used for the maximal
$\pi $-number dividing $n$. A subgroup $H$ of $G$ is called a {\it
$\pi$-Hall subgroup}, if
$\pi(H)\subseteq\pi$ and $\pi(|G:H|)\subseteq \pi'$. The set of all
$\pi$-Hall subgroups of~$G$ we denote by
$\operatorname{Hall}_\pi(G)$. A Hall subgroup is  $\pi$-Hall subgroup for
some~$\pi$.

According to  [7] we say that $G$ {\it
satisfies $E_\pi$} (or briefly  $G\in E_\pi$), if $G$ possesses a
$\pi$-Hall subgroup. If, moreover, every two $\pi$-Hall subgroups are conjugate,
then we say that  $G$ {\it satisfies
$C_\pi$} ($G\in C_\pi$). If, in addition, each $\pi$-subgroup of~$G$ is
included in a  $\pi$-Hall subgroup, then we say that $G$ {\it satisfies
$D_\pi$} ($G\in D_\pi$). A group satisfying
$E_\pi$ ($C_\pi$, $D_\pi$) we also call an
$E_\pi$- (respectively a $C_\pi$-, a $D_\pi$-) {\it group}.

A finite group possessing a (sub)normal series such that all factors of the
series are either $\pi$- or $\pi'$-groups is called $\pi$-{\it
separable}.

\begin{Lemma} {\rm [7, Lemma 1]}
Let $A$ be a normal subgroup of a finite group~$G$.  If
$G\in E_\pi$ and
$H\in\operatorname{Hall}_\pi(G)$, then
$A,G/A\in E_\pi$, moreover  $H\cap A\in \operatorname{Hall}_\pi(A)$ and
$HA/A\in\operatorname{Hall}_\pi(G/A)$.
\end{Lemma}

\begin{Lemma} {\rm [8; 7, Corollary D5.2]}
A $\pi$-separable group satisfies~$D_\pi$.
\end{Lemma}

\begin{Lemma}
Let $H$ be a subgroup of $G$, $g\in G$, $y\in \langle H,
H^g\rangle$. If subgroups $H^y$ and $H^g$ are conjugate in $\langle H^y,
H^g\rangle$, then $H$ and $H^g$ are conjugate in
$\langle H, H^g\rangle$.
\end{Lemma}

\begin{proof} Let $z\in \langle H^y, H^g\rangle$, and
$H^{yz}=H^g$. Then  $z\in \langle H, H^g\rangle$ since
$\langle H^y,H^g\rangle$
$\le\langle H, H^g\rangle $. Put $x=yz$.
Then $x\in\langle H,H^g\rangle$ and $H^x=H^g$.
\end{proof}

\begin{Lemma}
Let $H$ be a subgroup of a finite  $G$. Assume that  $H$ includes a pronormal
$($for example, a Sylow$)$  subgroup $S$ of
$G$. Then the following statements are equivalent:

$(1)$~$H\operatorname{prn} G;$

$(2)$~$H$ and $H^g$ are conjugate in $\langle H, H^g\rangle$ for each $g\in
N_G(S)$.
\end{Lemma}

\begin{proof}  Clearly $(1)\Rightarrow(2)$. We prove that
$(2)\Rightarrow(1)$. Assume that statement $(2)$ holds.
Choose arbitrary $g\in G$. Notice that $S,S^g\leq
\langle H, H^g\rangle$. Since $S$ is pronormal, there exists $y\in \langle S,
S^g\rangle\leq \langle H, H^g\rangle$ such that the equality $S^{gy}=S$ holds.
In particular,  $gy\in N_G(S)$. In view of $(2)$, subgroups
$H$ and $H^{gy}$ are conjugate in $\langle H, H^{gy}\rangle$. Then
$H^{y^{-1}}$ and $H^g$ are conjugate in $\langle H^{y^{-1}},
H^g\rangle$. Now $H$ and $H^g$ are conjugate in
$\langle H, H^g\rangle$ by Lemma~4.
\end{proof}

\begin{Lemma}
Let $\overline{\phantom{g}}:G\rightarrow G_1$ be a homomorphism of groups,
$H\le G$.   If $H\operatorname{prn} G$, then
$\overline{H}\operatorname{prn}\overline{G}$.
\end{Lemma}

\begin{proof} Clear.\end{proof}

\begin{Lemma}
Let $G$ be a finite group and $G_1,\dots,G_n$ be normal subgroups of $G$ such
that $[G_i,G_j]=1$ for
$i\ne j$ and $G=G_1\dots G_n$. Assume that for each
$i=1,\dots, n$ a pronormal subgroup $H_i$ of $G_i$ is chosen, and  $H=\langle
H_1,\dots, H_n\rangle$. Then
$H\operatorname{prn} G$.
\end{Lemma}

\begin{proof}
Choose arbitrary $g\in G$. Then $g=g_1\dots g_n$ for some
$g_1\in G_1, \dots, g_n\in G_n$. Since $H_i$ is pronormal in $G_i$ for each
$i=1,\dots, n$, there exist $x_i\in \big\langle H_i,H_i^{g_i}\big\rangle$ such
that $H_i^{x_i}=H_i^{g_i}$. Since  $[G_i,G_j]=1$ for
$i\ne j$, we have  $H_i^g=H_i^{g_i}$ for each $i=1,\dots ,n$. The same
arguments imply $H_i^{x_i}=H_i^{x}$, where
$x=x_1\dots x_n$. Clearly
$$
x\in \big\langle H_i,H_i^{g_i}\mid i=1,\dots,n\big\rangle
=\big\langle H_i,H_i^{g}\mid i=1,\dots,n\big\rangle=\langle H, H^g\rangle.
$$
Further,
\begin{multline}
H^g=\big\langle H_i^{g}\mid i=1,\dots, n\big\rangle
=\big\langle H_i^{g_i}\mid i=1,\dots,n\big\rangle
=\big\langle H_i^{x_i}\mid i=1,\dots, n\big\rangle
\\
=\big\langle H_i^{x}\mid i=1,\dots, n\big\rangle=H^x.
\end{multline}
\end{proof}

\begin{Lemma}
Let $G$ be a finite group, $H\in\operatorname{Hall}_\pi(G)$,
$A\trianglelefteq G$, and $G=HA$. If $(H\cap A)\operatorname{prn} A$, then
$H\operatorname{prn}G$.
\end{Lemma}

\begin{proof}
By Lemma~2, $H\cap A$ is a $\pi$-Hall subgroup of $A$. Let $(H\cap
A)\operatorname{prn} A$. Choose arbitrary
$g\in G$ and show that $H^x=H^g$ for some
$x\in\langle H, H^g\rangle$.

Since $G=HA$, there exist $h\in H$ and $a\in A$ such that $g=ha$. Since
$(H\cap A)\operatorname{prn} A$, there exists $y\in   \langle
H\cap A, H^a\cap A\rangle$ such that $H^y\cap A=H^a\cap A$.
In view of
$$
y\in \langle H\cap A, H^a\cap A\rangle\le\langle H, H^a\rangle
=\langle H,H^{ha}\rangle=\langle H, H^g\rangle,
$$
and Lemma~4 we need to consider the case $H=H^y$. In particular,
$$
H\cap A=H^a\cap A=H^g\cap A.
$$
Now
$H$, $H^g$, and $g$ are included in $N_G(H\cap A)$. Since $G=HA$ we have
$G=AN_G(H\cap A)$. Notice that
$$
N_G(H\cap A)/N_A(H\cap A)=N_G(H\cap A)/(A\cap N_G(H\cap A))
\simeq AN_G(H\cap A)/A=G/A
$$
is a $\pi$-group.
Consider a normal series
$$
N_G(H\cap A)\trianglerighteq  N_A(H\cap A)
\trianglerighteq  H \cap A\trianglerighteq  1
$$
of $N_G(H\cap A)$. Each factor of the series is either a $\pi$- or a
$\pi'$-group, so $N_G(H\cap A)$ is  $\pi$-separable. Therefore,  the subgroup
$\langle H,H^g\rangle$ of $N_G(H\cap A)$ is
$\pi$-separable as well, and in particular
$\langle H,H^g\rangle\in D_\pi$ by Lemma~3. Thus
$\pi$-Hall subgroups $H$ and $H^g$ are conjugate in $\langle
H,H^g\rangle$.
\end{proof}

The next lemma gives a sufficient condition for the treatment of
lemma~6 in case when $H$ is a Hall subgroup of~$G$.

\begin{Lemma}
Let $\frak{X}$ be a class of finite groups close under subgroups such
that $\frak{X}\subseteq C_\pi$. Let $G$ be a finite group,
$H\in\operatorname{Hall}_\pi(G)$, $A\trianglelefteq G$, and
$\overline{\phantom{g}}:G\rightarrow G/A$ be the natural homomorphism.
Assume also that $A\in\frak{X}$. Then $H\operatorname{prn} G$ if and only if
$\overline{H}\operatorname{prn} \overline{G}$.
\end{Lemma}

\begin{proof} The implication $\Rightarrow$ holds by Lemma~6.

We prove $\Leftarrow$. Let $g\in G$. We need to show that $H^x=H^g$ for some
$x\in\langle H, H^g\rangle$.  Since $\overline{H}\operatorname{prn}
\overline{G}$, there exists $y\in\langle H, H^g\rangle$ such that
$H^yA=H^gA$. By Lemma~4 we may substitute $H$ by $H^y$ and so we may assume that
$HA=H^gA$.

Consider $M=\langle H\cap A, H^g\cap A\rangle$.  Since
$M \le A$, $A\in \frak{X}$  and $\frak{X}$ is closed under subgroups, we have
\hbox{$M\in\frak{X}\subseteq C_\pi$.} Further
$H\cap A, H^g\cap A\in\operatorname{Hall}_\pi(A)$ by Lemma~2, and
$M\leq A$, so  $H\cap A, H^g\cap
A\in\operatorname{Hall}_\pi(M)$. Hence $H^a\cap A= H^g\cap A$ for some
$a\in M$. Since $M\le \langle H, H^g\rangle$, by Lemma~4  we may substitute
$H$ by $H^a$, and so we may assume that $H\cap A=H^g\cap A$. In such case
$g\in N_G(H\cap A)$ and $H,H^g\le N_G(H\cap A)$.  Since
$A\in C_\pi$ by Frattini argument we have
$G=AN_G(H\cap A)$. Now
$$
N_G(H\cap A)/N_A(H\cap A)=N_G(H\cap A)/(A\cap N_G(H\cap A))
\simeq AN_G(H\cap A)/A=G/A=\overline{G}.
$$
As we noted above
$\overline{H}=\overline{H}^{\overline{g}}$, so the isomorphism implies that
$HN_A(H\cap A)=H^gN_A(H\cap A)$. Denote the last subgroup by $B$ for brevity.
Then $B$ is $\pi$-separable and $H,H^g\le B$. Moreover,
$\langle H, H^g\rangle$ is also $\pi$-separable as a subgroup of a
$\pi$-separable group~$B$. In particular, by Lemma~3
$$
\langle H, H^g\rangle\in D_\pi\quad\text{ and }\quad
H,H^g\in\operatorname{Hall}_\pi(\langle H, H^g\rangle),
$$
whence $H$ and $H^g$ are conjugate in $\langle H, H^g\rangle$.
\end{proof}

Let $G$ be a finite group and 1$\pi(G)=\{p_1,\dots, p_n\}$. Following
[7] we say that $G$ has a {\it Sylow series of
complexion}\footnote{Parentheses in the notation
$(p_1,\dots, p_n)$ are used for an ordered set, apart from braces. For
example, the symmetric group $\operatorname{Sym}_3$ has a Sylow series of
complexity $(2,3)$,  while the alternating group $\operatorname{Alt}_4$ has a
Sylow series of complexity  $(3,2)$.}
$(p_1,\dots, p_n)$, if $G$ possesses a normal series
$$
G=G_0>G_1>\dots>G_n=1
$$
such that each section $G_{i-1}/G_i$ is isomorphic to a Sylow
$p_i$-subgroup of~$G$.

\begin{Lemma}
Let $G$ be a finite group, $H$ be its Hall subgroup with a Sylow series. Then
$H\operatorname{prn} G$.
\end{Lemma}

\begin{proof} Let $g\in G$. We show that $H$ and $H^g$ are conjugate in
$\langle H,H^g\rangle$. By  [7, Theorem A1] every two Hall subgroups of a
finite groups having Sylow series of the same complexion are conjugate. Since
$H$ and $H^g$ are two Hall subgroups of
$\langle H,H^g\rangle$ having Sylow series of the same complexion,  $H$ and
$H^g$ are conjugate in $\langle H,H^g\rangle$.
\end{proof}

\begin{Lemma}
Let  $G$ be a finite nonabelian simple group, $H$ be its Hall subgroup of order
not divisible either by $2$ or by $3$.
Then $H$ has a Sylow series.
\end{Lemma}

\begin{proof} If $2$ does not divide the order of $H$, the claim is proven in
[9, Theorem~B]. If $3$ does not divide the order of$H$ the claim follows from
[10, Lemma~5.1, Theorem~5.2].
\end{proof}

The symmetric group and the alternating group of
degree $n$ we denote by $\operatorname{Sym}_n$ and $\operatorname{Alt}_n$
respectively.

A finite field containing $q$ elements, is denoted by~${\Bbb F}_q$.

Given odd number $q$ define $\varepsilon(q)=(-1)^{(q-1)/2}$, i.~e.,
$\varepsilon(q)=1$, if $q-1$ is divisible b~$4$, and $\varepsilon(q)=-1$
otherwise. Without additional explanations we use symbols $\varepsilon,\eta$ to
denote either an element from $\{+1,-1\}$ or the sign of the element.

Given group of Lie type the order of the base field is always denoted by $q$
(see [1], for example), while its characteristic is denoted by $p$.
Given matrix group $G$ the reduction modulo scalars is denoted by
${\mathrm{P}} G$.

Our notation for classical groups agrees with that of [11]. We recall special
notation, that we often use:

 $\operatorname{GL}^+_n(q)=\operatorname{GL}_n(q)$ is a
  general linear group of degree $n$ over~${\Bbb F}_q$;

$\operatorname{SL}^+_n(q)=\operatorname{SL}_n(q)$ is a
  special linear group of degree $n$ over~${\Bbb F}_q$;

$\operatorname{PGL}^+_n(q)=\operatorname{PGL}_n(q)$ is a
  projective general linear group of degree $n$ over~${\Bbb F}_q$;

$\operatorname{PSL}^+_n(q)=\operatorname{PSL}_n(q)$ is a
  projective special linear group of degree~$n$ over~${\Bbb F}_q$;

$\operatorname{GL}^-_n(q)=\operatorname{GU}_n(q)$ is a
  general unitary group of degree $n$ over~${\Bbb F}_{q^2}$;

$\operatorname{SL}^-_n(q)=\operatorname{SU}_n(q)$ is a
  special unitary group of degree $n$ over~${\Bbb
F}_{q^2}$;

$\operatorname{PSL}^-_n(q)=\operatorname{PSU}_n(q)$ is a
  projective special unitary group of degree  $n$ over~${\Bbb F}_{q^2}$;

$\operatorname{PGL}^-_n(q)=\operatorname{PGU}_n(q)$ is a
  projective general unitary group of degree $n$ over~${\Bbb F}_{q^2}$;

$\operatorname{Sp}_n(q)$ is a
  simplectic group of degree  $n$ over~${\Bbb F}_q$;

 $\operatorname{PSp}_n(q)$ is a
  projective simplectic group of degree $n$ over~${\Bbb F}_q$.

Necessary facts about properties and  structure of finite groups of
Lie type can be found in [12--15], properties and structure of linear algebraic
groups can be found in~[12], results concerning the connection between groups
of Lie type and linear algebraic groups can be found in [13--14]. Also in
[13--14] the definitions of Borel an Cartan subgroups, a parabolic subgroup, and
a maximal torus in a finite group of Lie type can be found.

We denote groups $E_6(q)$ and ${{}^2E}_6(q)$ by
$E^+_6(q)$ and $E^-_6(q)$ respectively.

A {\it Frobenius map} of an algebraic group $\overline{G}$ is a surjective
endomorphism  $\sigma:\overline{G}\rightarrow \overline{G}$ such that the
set of its stable points $\overline{G}_\sigma$ is finite. Each simple
group of Lie type of a finite field $F$ of characteristic $p$ is known to
coincide with $O^{p'}(\overline{G}_\sigma)$ for an appropriate linear algebraic
group $\overline{G}$ over the algebraci closure of $F$ and a Frobenius
map~$\sigma$, where $O^{p'}(\overline{G}_\sigma)$ is a  subgroup of
$\overline{G}_\sigma$ generated by all $p$-elements.

Let $\overline{R}$ be a closed $\sigma$-stable subgroup of an algebraic group
$\overline{G}$ for a Frobenius map $\sigma$ of
$\overline{G}$. Consider subgroups  $R=G\cap \overline{R}$ and $N(G,R)=G\cap
N_{\overline{G}}(\overline{R})$, where $G=O^{p'}(\overline{G}_\sigma)$. Notice
that  $N(G,R)\le N_G(R)$ and $N(G,R)\not=N_G(R)$ in general.

\begin{Lemma} {\rm [16, Corollary of Theorems 1--3]}
Let $G$ be a finite simple nonabelian group and
$S\in\operatorname{Syl}_2(G)$. Then $N_G(S)=S$, except the following cases:

$(1)$~$G\simeq  J_2,~J_3,~Suz$ or $HN$ and $|N_G(S):S|=3$;

$(2)$~$G\simeq  {^2}G_2(q)$ or $J_1$ and
$N_G(S)\simeq 2^3.7.3<\operatorname{Hol}(2^3)$;

$(3)$~$G$ is a group of Lie type over a field of characteristic $2$ and
$N_G(S)$ is a Borel subgroup of~$G$;

$(4)$~$G\simeq \operatorname{PSL}_2(q)$, where $3<q \equiv  \pm 3\pmod 8$ and
$N_G(S)\simeq \operatorname{Alt}_4$;

$(5)$~$G \simeq E^{\eta}_6(q)$, $\eta=\pm$, $q$ is odd and
$N_G(S)=S\times C$, where $C$ is a nontrivial cyclic group of order
$(q-\eta)_{2'}/(q-\eta,3)_{2'}$;

$(6)$~$G\simeq\operatorname{PSp}_{2m}(q)$, $m\geq 2$, $q\equiv\pm 3\pmod 8$,
the factor group $N_G(S)/S$ is isomorphic to an elementary abelian  $3$-group
of order $3^t$ and $t$ can be found from the $2$-adic decomposition
$$
m=2^{s_1}+\dots+2^{s_t},
$$
where $s_1>\dots>s_t\geq 0$;

$(7)$~$G\simeq \operatorname{PSL}_n^\eta(q)$, $n\geq 3$, $\eta=\pm$, $q$ is
odd,
$$
N_G(S)\simeq S\times C_1\times \dots\times C_{t-1},
$$
$t$ can be found from a $2$-adic decomposition
$$
n=2^{s_1}+\dots+2^{s_t},
$$
where $s_1>\dots>s_t\ge0, $ and $C_1,\dots,C_{t-2},C_{t-1}$ are cyclic groups
of orders
${(q-\eta)}_{2'},\ldots,{(q-\eta)}_{2'}$, ${(q-\eta)}_{2'}/{(q-\eta,n)} _{2'}$
respectively.
\end{Lemma}

\begin{Lemma} {\rm [7, Theorem~A4; 17]}
Let $2,3\in\pi$. Then the list of all cases, when  $\operatorname{Sym}_n$
possesses a proper $\pi$-Hall subgroup is given in Table~$1$. In particular,
each proper $\pi$-Hall subgroup of $\operatorname{Sym}_n$ is maximal in
$\operatorname{Sym}_n$.
\end{Lemma}

\begin{longtable}{|c|c|c|}\caption{$\pi$-Hall subgroups in symmetric groups\label{Symmetric}}\\  \hline
$n$&$\pi\cap\pi(\Sym_n)$&$H\in\Hall_\pi(\Sym_n)$\\ \hline\hline
prime&$\pi((n-1)!)$&$\Sym_{n-1}$\\
$7$&$\{2,3\}$&$\Sym_3\times\Sym_4$\\
$8$&$\{2,3\}$&$\Sym_4\wr\Sym_2$\\ \hline
\end{longtable}

\begin{Lemma} {\rm [18, Theorem 4.1]}
Let $G$ be either one of~$26$ sporadic groups, or the Tits group. Assume that
$\pi$ contains both $2$ and $3$. Then  $G$ possesses a proper
$\pi$-Hall subgroup~$H$ if and only if one of the conditions on  $G$ and
$\pi\cap\pi(G)$ from Table.~$2$ holds. In the table the structure of  $H$ is
also given.
\end{Lemma}

\begin{longtable}{|l|l|r|}\caption{$\pi$-Hall subgroups in sporadic groups,
case $2,3\in\pi$}\label{tb0}\\  \hline
$G$ & $\pi\cap\pi(G)$ &
Structure $H$ \\ \hline\hline $M_{11}$  & $\{2,3\}$ & $3^2\splitext
Q_8\arbitraryext 2$ \\ &
$\{2,3,5\}$       &  $\Alt_6\arbitraryext 2$\\ \hline $M_{22}$  & $\{2,3,5\}$ &
$2^4\splitext \Alt_6$\\ \hline $M_{23}$  & $\{2,3\}$         &
$2^4\splitext(3\times
A_4)\splitext2$\\ & $\{2,3,5\}$ &  $2^4\splitext \Alt_6$\\ & $\{2,3,5\}$       &
$2^4\splitext (3\times \Alt_5)\splitext 2$\\ & $\{2,3,5,7\}$     &  ${\rm
          L}_3(4)\splitext 2_2$\\ & $\{2,3,5,7\}$     &  $2^4\splitext \Alt_7$\\
&
          $\{2,3,5,7,11\}$  & $M_{22}$\\ \hline $M_{24}$  & $\{2,3,5\}$
          & $2^6 \splitext 3\nonsplitext \Sym_6$\\ \hline $J_1$     & $\{2,3\}$
       &
$2\times \Alt_4$\\ & $\{2,3,5\}$ &
          $2\times A_5$\\ & $\{2,3,7\}$       &  $2^3\splitext 7\splitext 3$\\
\hline
$J_4$     & $\{2,3,5\}$       &  $2^{11}\splitext (2^6\splitext 3\nonsplitext
\Sym_6)$\\
\hline
\end{longtable}

\begin{Lemma} {\rm [19, Theorem~3.3]}
Let $G$ be a finite group of Lie type over a field of characteristic  $p\in\pi$.
If $H$ is a $\pi$-Hall subgroup of $G$, then either $H$ is included in a Borel
subgroup, or $H$ is a parabolic subgroup of~$G$.
\end{Lemma}

\begin{Lemma} {\rm [20, Lemma 3.1]}
Let  $G\simeq{\mathrm{P}}\operatorname{SL}_2(q)
\simeq{\mathrm{P}}\operatorname{SL}_2^\eta(q)
\simeq{\mathrm{P}}\operatorname{Sp}_2(q)$, where
$q$ is a power of an odd prime $p$, and set
$\varepsilon=\varepsilon(q)$. Assume that $2,3\in\pi$, and
$p\not\in\pi$.
Then $G\in E_\pi$ if and only of on of the cases from Table~$3$ holds.
\end{Lemma}

\begin{longtable}{|c|c|c|}\caption{$\pi$-Hall subgroups $H$ of
$\PSL_2(q)$, $2,3\in\pi$, $p\not\in\pi$}\label{SL23Halldim2}\\  \hline
$\pi\cap \pi(G)$& $H$&conditions \\ \hline\hline
$\subseteq\pi(q-\varepsilon)$&$D_{q-\varepsilon}$&---\\ \hline
$\{2,3\}$&$\Alt_4$&$(q^2-1)_{\mbox{}_{\{2,3\}}}=24$\\ \hline
$\{2,3\}$&$\Sym_4$&$(q^2-1)_{\mbox{}_{\{2,3\}}}=48$\\ \hline
$\{2,3,5\}$&$\Alt_5$&$(q^2-1)_{\mbox{}_{\{2,3,5\}}}=120$\\ \hline
\end{longtable}

\begin{Lemma} {\rm[20, Lemma 3.2]}
Assume that $G=\operatorname{GL}^\eta_2(q)$, where $q$ is a power of a prime
$p$, ${\mathrm{P}}: G\rightarrow G/Z(G) =\operatorname{PGL}_2^\eta(q)$ is the
natural homomorphism, and let  $\varepsilon=\varepsilon(q)$. Assume also that
$2,3\in\pi$ and $p\not\in\pi$. A subgroup $H$ of $G$ is a $\pi$-Hall
subgroup if and only if one of the following statements holds:

$(1)$~$\pi\cap\pi(G)\subseteq \pi(q-\varepsilon)$,
${\mathrm{P}} H$ is a $\pi$-Hall subgroup of the dihedral group
$D_{2(q-\varepsilon)}$ of order
$2(q-\varepsilon)$ of~${\mathrm{P}} G$;

$(2)$~$\pi\cap\pi(G)=\{2,3\}$, $(q^2-1)_{{}_{\{2,3\}}}=24$,
${\mathrm{P}} H\simeq \operatorname{Sym}_4$.
\noindent

Moreover every two $\pi$-Hall subgroups of $G$, satisfying to the same
statement  $(1)$ or $(2)$ are conjugate.
\end{Lemma}

\begin{Lemma} {\rm [20, Lemma~4.3]}
Let $G^*=\operatorname{SL}_n^\eta(q)$ be a special linear or unitary
group with the base field ${\Bbb F}_q$ of characteristic $p$, and let $n\ge2$.
Assume that
$2,3\in\pi$ and $p\not\in\pi$. Suppose that
$G^*\in E_\pi$ and $H^*$ is a $\pi$-Hall subgroup of $G^*$.
Then for $G^*$, $H^*$ and $\pi$ one of the following statements holds.

$(1)$~$n=2$ and for groups $G=G^*/Z(G^*)$ and $H=H^*Z(G^*)/Z(G^*)$
the conditions from Table~$3$ holds.

$(2)$~Either $q\equiv \eta \pmod {12}$, or $n=3$ and
 $q\equiv \eta \pmod 4$;  $\operatorname{Sym}_n$ satisfies $E_\pi$,
$\pi\cap\pi(G^*)\subseteq \pi(q-\eta)\cup \pi(n!)$ and if
$r\in(\pi\cap\pi(n!))\setminus
\pi(q-\eta)$, then $|G^*|_r=|\operatorname{Sym}_n|_r$;
$H^*$ is included in
$$
M=L\cap G^*\simeq Z^{n-1}\, .\, \operatorname{Sym}_n,
$$
where $L=Z\wr \operatorname{Sym}_n\le \operatorname{GL}_n^\eta(q)$ and
$Z=\operatorname{GL}_1^\eta(q)$ is a cyclic group of order
$q-\eta$.

$(3)$~$n=2m+k$, where $k\in\{0,1\}$, $m\geq 1$,
$q\equiv -\eta \pmod {3}$, $\pi\cap\pi(G^*)\subseteq \pi(q^2-1)$, both
$\operatorname{Sym}_m$
and $\operatorname{GL}_2^\eta(q)$ satisfy $E_\pi$
\footnote{By Lemma~16 conditions
$\operatorname{GL}_2^\eta(q)\in E_\pi$ and $q\equiv -\eta\pmod3$
mean that $q\equiv -\eta\pmod{r}$ for all odd primes
$r\in\pi(q^2-1)\cap\pi$.};
$H^*$ is included in
$$
M=L\cap G^*\simeq\big(\underbrace{\operatorname{GL}_2^\eta(q)
\circ\dots\circ\operatorname{GL}_2^\eta(q)}_{
m \text{ \rm times}}\big)\, .\,\operatorname{Sym}_m\circ Z,
$$
where $L=\operatorname{GL}_2^\eta(q)\wr
\operatorname{Sym}_m\times Z\leq \operatorname{GL}_n(q)$
and $Z$ is a cyclic group of order
${q-\eta}$ if  $k=1$, and  $Z=1$ if $k=0$.
A subgroup  $H^*$ acting by conjugation on the set of factors of type
$\operatorname{GL}_2^\eta(q)$ in the central product
$$
\underbrace{\operatorname{GL}_2^\eta(q)\circ\dots\circ\operatorname{GL}_2^\eta(q)}_{
m \text{ \rm times}},
$$
has at most two orbits. The intersection of  $H^*$ with each factor
$\operatorname{GL}_2^\eta(q)$ in
$(1)$ is a $\pi$-Hall subgroup of
$\operatorname{GL}_2^\eta(q)$. All intersections of $H^*$ with factors from the
same orbit satisfy to the same statement $(1)$ or $(2)$ in Lemma~$17$.

$(4)$~$n=4$, $\pi\cap \pi(G^*)=\{2,3,5\}$,  $q\equiv5\eta\pmod 8$,
$(q+\eta)_3=3$,
$(q^2+1)_5=5$ and $H^*\simeq 4\, .\, 2^4\, .\, \operatorname{Alt}_6$.

$(5)$~$n=11$, $\pi\cap \pi(G^*)=\{2,3\}$,
$(q^2-1)_{\text{}_{\{2,3\}}}=24$, $q\equiv -\eta\pmod3$, $q\equiv\eta\pmod4$,
$H^*$ is included in $M=L\cap G^*,$ where $L$ is a subgroup of
$G^*$ of type $\big(\big(\operatorname{GL}_2^\eta(q)\wr
\operatorname{Sym}_4\big)
\perp \big(\operatorname{GL}_1^\eta(q)\wr
\operatorname{Sym}_3\big)\big)$
and
$$
H^*= (((Z\circ 2\, .\, \operatorname{Sym}_4)
\wr \operatorname{Sym}_4)\times (Z\wr \operatorname{Sym}_3))\cap G,
$$
where $Z$ is a Sylow $2$-subgroup of a cyclic group of order~${q-\eta}$.
\end{Lemma}

\begin{Lemma} {\rm [20, Lemma 4.4]}
Let $G^*=\operatorname{Sp}_{2n}(q)$ be a simplectic group over a field
${\Bbb F}_q$ of characteristic $p$. Assume that $2,3\in\pi$ and
$p\not\in\pi$. Suppose that $G^*\in E_\pi$ and
$H^*\in \operatorname{Hall}_\pi(G)$.
Then both $\operatorname{Sym}_n$ and $\operatorname{SL}_2(q)$ satisfy~$E_\pi$
and $\pi \cap\pi(G^*)\subseteq \pi(q^2-1)$. Moreover $H^*$ is a
$\pi$-Hall subgroup of
$$
M=\operatorname{Sp}_2(q)\wr \operatorname{Sym}_n\simeq
(\underbrace{\operatorname{SL}_2(q)\times\dots
\times\operatorname{SL}_2(q)}_{n\text{ \rm times}})
:\operatorname{Sym}_n\leq G^*.
$$
\end{Lemma}

\begin{Lemma} {\rm [20, Lemma~7.3]}
Let $G=E_6^\eta(q)$, where $q$ is a power of a prime $p$, and
$\varepsilon=\varepsilon(q)$. Assume that $2,3\in\pi$ and
$p\not\in\pi$. Suppose that $G$ possesses a
$\pi$-Hall subgroup $H$. Then $\pi\cap \pi(G)\subseteq \pi(q-\varepsilon)$ and
one of the following statements holds:

$(1)$~$\eta=\varepsilon$, $5\in\pi$ and $H$ is a $\pi$-Hall subgroup of
$$
M={((q-\eta)^6.W(E_6))}/(3,q-\eta);
$$

$(2)$~$\eta=-\varepsilon$ and $H$ is a $\pi$-Hall subgroup of
$$
M={(q^2-1)^2(q+\eta)}^2.W(F_4).
$$
\end{Lemma}

\section{Proof of Theorem~1}

Let $G$ be a finite simple group and
$H\in\operatorname{Hall}_\pi(G)$. We show that
$H\operatorname{prn} G$, and thus we prove Theorem~1. By Lemmas~10 and~11 we
may assume that $2,3\in\pi$. Let
$S\in\operatorname{Syl}_2(H)\subseteq\operatorname{Syl}_2(G)$
and $g\in N_G(S)$ be arbitrary. By Lemma~5 it is enough to prove that
$H$ and $H^g$ are conjugate in $\langle H,H^g\rangle$. If
$N_G(S)=S$, then this statement is true:
$g\in N_G(S)=S\le H$, so $H^g=H$. Therefore we may assume that one of the
exceptional cases (1)--(7) from Lemma~12 holds, and $H$ is a proper
$\pi$-Hall subgroup of~$G$.

We consider cases (1)--(7) from Lemma~12, proving a series of auxiliary lemmas.
In order to unify the notation in Lemmas with already introduced notation we
say that $(\star)$ {\it holds}, if

(a)~$G$ is a finite simple group;

(b)~$2,3\in\pi$;

(c)~$H\in\operatorname{Hall}_\pi(G)$ and $H<G$;

(d)~$S\in\operatorname{Syl}_2(H)\subseteq\operatorname{Syl}_2(G)$;

(e)~$g\in N_G(S).$

The following lemma follows from Lemma~14 immediately.

\begin{Lemma}
Assume that $(\star)$ holds. If
$G\simeq  J_2, J_3, Suz$ or $HN$, then
$G$ does not possesses proper $\pi$-Hall subgroups.
\end{Lemma}

Thus if case~(1) of Lemma~12 holds, then by Lemma~5 $H\operatorname{prn} G$.

\begin{Lemma} Assume that $(\star)$ holds. Then the following statements
hold.

$(1)$~If $G\simeq {}^2G_2(q)$, then $G$ does not possesses proper
$\pi$-Hall subgroup.

$(2)$~If $G\simeq J_1$, then one of the following cases holds:

$\quad${\rm (a)}~$H\simeq 2\times\operatorname{Alt}_4$ and $H$ possesses a
Sylow tower;

$\quad${\rm (b)}~$H\simeq 2^3:7:3$ and $H$ possesses a Sylow tower;

$\quad${\rm (c)}~$H\simeq 2\times \operatorname{Alt}_5$ and $H$ is maximal
in~$G$.

$(3)$~If  $G\simeq J_1$, then $H$ is conjugate with $H^g$ by an element from
$\langle H,H^g\rangle$.
\end{Lemma}

\begin{proof}
Statement  (1) follows from [19, Theorem~1.2], since
$3\in\pi$ and $3$ is the characteristic of the base field for
${}^2G_2(q)$. Lemma~14 implies the structure of $H$ in statement (2), moreover
it is clear that in cases  (a) and (b) the subgroup has a Sylow series, In case
(c) $H$ is maximal in view of [6]. Statement (3) follows from (2), Lemma~10,
and pronormality of maximal subgroups.
\end{proof}

Thus Lemmas~5 and~22 imply that
$H\operatorname{prn} G$, if statement~(2) of Lemma~12 holds.

\begin{Lemma}
Assume that $(\star)$ holds, and $G$ is a group of Lie type over a field of
characteristic $2$.  Then $S$ is a maximal unipotent subgroup,
$N_G(S)$ is a Borel subgroup of $G$ and one of the following statements holds:

$(1)$~$H$ is included in a Borel subgroup and has a Sylow series;

$(2)$~$H$ is parabolic and includes $N_G(S)$.

In both cases $H$ is conjugate with $H^g$ by an element from $\langle
H,H^g\rangle$.
\end{Lemma}

\begin{proof}
In view of Lemma~15, the structure of Borel subgroups and the fact that every
parabolic subgroup includes a Borel subgroup we obtain that either (1) or (2)
holds. By using Lemma~10 we conclude that  $H$ is conjugate with $H^g$ by an
element from $\langle H,H^g\rangle$, if statement $(1)$ holds. If statement
$(2)$ holds the final conclusion is evident, since $g\in H$.
\end{proof}

Thus if statement~(3) of Lemma~12 holds, then
$H\operatorname{prn} G$.

\begin{Lemma}
Assume that $(\star)$ holds,
$G=\operatorname{PSL}_2(q)$, $q\equiv \pm 3\pmod8$, and $q>3$.
Then one of the following statements holds:

$(1)$~$H$ is a $\pi$-Hall subgroup in a dihedral group of order
$q-\varepsilon$, where $\varepsilon=\varepsilon(q)=(-1)^{(q-1)/2}$, and it has
a Sylow series;

$(2)$~$H\simeq \operatorname{Alt}_4$ and $H$ has a Sylow sereis;

$(3)$~$H\simeq  \operatorname{Alt}_5$ and $H$ includes
$N_G(S)\simeq \operatorname{Alt}_4$, in particular, $H^g=H$.

In any case $H$ is conjugate with  $H^g$ by an element from $\langle
H,H^g\rangle$.
\end{Lemma}

\begin{proof}
Conditions $q\equiv \pm 3\pmod8$ and $q>3$, and Lemma~16 imply the
structure of $H$. Moreover, if either $H$ is included in a dihedral subgroup,
or   $H\simeq \operatorname{Alt}_4$,  then it clearly has a Sylow series.
Assume that $H\simeq \operatorname{Alt}_5$.  Then
$\operatorname{Alt}_4=N_H(S)\le N_G(S)\simeq \operatorname{Alt}_4$,
so $N_H(S)=N_G(S)$.  Using statements
$(1)$--$(3)$ and Lemma~10 we obtain the final conclusion.
\end{proof}

Thus we have shown that if statement~(4) of Lemma~12 holds, then
$H\operatorname{prn} G$.

Lemma~24 implies also the following statement that is extensively used for
consideration of items (6) and (7) in Lemma~12.

\begin{Lemma}
Let $2,3\in\pi$,  $q$ be a power of an odd prime $p\not\in\pi$,
$$G^*\in\big\{\operatorname{PSL}_2(q),\ \operatorname{PGL}^\eta_2(q),\
\operatorname{SL}_2(q),\ \operatorname{GL}^\eta_2(q)\big\}$$
and $H^*\in\operatorname{Hall}_\pi(G^*)$.
Then $H^*\operatorname{prn} G^*$.
 \end{Lemma}

\begin{proof}
If $G^*=\operatorname{PSL}_2(q)$
and $S^*\in \operatorname{Syl}_2(H^*)\subseteq \operatorname{Syl}_2(G^*)$,
then by Lemma~12 either
$N_{G^*}(S^*)=S^*$ or $G^*$ satisfies the conditions of Lemma~24. In both cases
$H^*$ is pronormal.

Now let $G^*=\operatorname{PGL}_2^\eta(q)$ and
$A^*=\operatorname{PSL}^\eta_2(q)\simeq \operatorname{PSL}_2(q)$ be a normal
subgroup of index  $2$ in $G^*$.  As we have already shown, $H^*\cap
A^*\operatorname{prn} A^*$ and $G^*=A^*H^*$. Using Lemma~8 we conclude that
$H^*\operatorname{prn} G^*$.

Assume finally that $G^*$ is isomorphic to either
$\operatorname{SL}_2(q)$ or $\operatorname{GL}^\eta_2(q)$. Choose in Lemma~9
the class of all $2$-groups as ${\frak X}$. Then this lemma and the above
arguments imply $H^*\operatorname{prn} G^*$.
\end{proof}

Consider statement (5) in Lemma~12.

\begin{Lemma}
Assume that $(\star)$ holds and
$G=E_6^\eta(q)$, where  $q$ is a prime of
$p\not\in\pi$. Denote $\varepsilon(q)$ by $\varepsilon$. Then

$(1)$~$G$ includes an $S$-invariant maximal torus $T$ such that
  $$|T|=\left\{
  \begin{array}{ll}
  (q-\varepsilon)^6/(3,q-\varepsilon), & \text{\,если\,\,} \eta=\varepsilon,\\
  (q-\varepsilon)^4(q+\varepsilon)^2, & \text{\,если\,\,} \eta=-\varepsilon,\\
  \end{array}
  \right.
  $$ 
moreover $H\leq N_G(T)$ and $N_G(T)$ is an extension of $T$ by a $\pi$-group;

$(2)$~$N_G(T)$ includes $N_G(S)$;

$(3)$~subgroups $H$ and $H^g$ are conjugate in $\langle H,H^g\rangle$.
\end{Lemma}

\begin{proof}
(1)~The existence of $S$-invariant torus $T$ follows from [12, Theorem~4.10.2].
In view of [10, Lemma~3.10] such torus is unique up to conjugation and
$N_G(T)=N(G,T)$. Moreover by
[10, Lemma~3.11] the order of $T$ equals
$(q-\varepsilon)^6/(3,q-\varepsilon)$, if $\varepsilon=\eta$, and it equals
$(q-\varepsilon)^4(q+\varepsilon)^2$, if $\varepsilon=-\eta$.
Since $G\in E_\pi$ and $2,3\in\pi$, while $p\not\in\pi$, by Lemma~20 we obtain
that $H$ lies in $N_G(T)$ for some such torus $T$ and $N_G(T)/T$ is a
$\pi$-group.

(2)~Since $H\leq N_G(T)$ and $3\in\pi$, $N_G(T)$ includes a Sylow $3$-subgroup
of $G$. So in follows by [10, Lemma~3.13] that $N_G(S)\leq N_G(T)$.

(3)~In view of statement (2) of the lemma we remain to prove that $H
\operatorname{prn} N_G(T)$.  By~(1), $N_G(T)$ is an extension of an abelian
group $T$ by a $\pi$-group, in particular $N_G(T)=HT$.  Now, by Lemma~8,
${H \operatorname{prn} N_G(T)}$.
\end{proof}

Therefore, if statement~(5) of Lemma~12 holds, then
$H\operatorname{prn} G$.

In the next lemma we consider statement (6) and, partially, statement (7) of
Lemma~12.  We need to recall the notion of a fundamental subgroup introduced in
[21]. We use the notion in simple linear, unitary, and simplectic groups
in odd characteristic only, and their central extensions. Recall that
if $G$ is one of such groups, $X^+$ is a long root subgroup of~$G$, and
$X^-$ is the opposite root subgroup, then every $G$-conjugate of $\langle
X^+,X^-\rangle \simeq\operatorname{SL}_2(q)$ is called a {\it
fundamental subgroup}. If $S\in\operatorname{Syl}_2(G)$, then by
$\operatorname{Fun}_G(S)$ the sen of all fundamental subgroups $K$ of
$G$ such that  $K\cap S\in\operatorname{Syl}_2(K)$ is denoted.
$\operatorname{Fun}_G(S)$ is known  to be a maximal by inclusion
$S$-invariant set of pairwise commuting fundamental subgroups of  $G$ (see
[21]).

\begin{Lemma}
Assume that $(\star)$ holds and  $G$ is isomorphic to either
$\operatorname{PSL}_n^\eta(q)$ or
$\operatorname{PSp}_{n}(q)$, where $n>2$.
Let $\Delta=\operatorname{Fun}_G(S)$ and suppose that
$\Delta$ is $H$-invariant $($i.~e. $H\leq N_G(\Delta)$
in the notations of $[21])$.  Then  $H$ and  $H^g$ are conjugate in
$\langle H,H^g\rangle$.
 \end{Lemma}

 \begin{proof} Let $m=[n/2]$. Then $|\Delta|=m$.

In view of [11, Propositions 4.1.4, 4.2.9, and 4.2.10] the stabilizer
$N_G(\Delta)$ in $G$ of $\Delta$ coincides with the image in  $G$ of
a subgroup $M$ of either  $\operatorname{SL}_n^\eta(q)$ or
$\operatorname{Sp}_n(q)$, where $M$ is defined in the following way.
If $G=\operatorname{PSL}_n^\eta(q)$, then
$$
M= L\cap \operatorname{SL}_n^\eta(q)
\simeq\bigl(\underbrace{\operatorname{GL}_2^\eta(q)\circ
\dots\circ\operatorname{GL}_2^\eta(q)}_{m \text{ \rm times}}\bigr)\, .\,
\operatorname{Sym}_m\circ Z,
$$
moreover, $L=\operatorname{GL}_2^\eta(q)\wr \operatorname{Sym}_m
\times Z\le \operatorname{GL}_n^\eta(q)$, and $Z$ is a cyclic group of
order $q-\eta$ if $n$ is odd, and $Z=1$ if $n$ is even.  If
$G=\operatorname{PSp}_n(q)$, then
$$
M=\operatorname{Sp}_2(q)\wr \operatorname{Sym}_m\simeq
(\underbrace{\operatorname{SL}_2(q)\times\dots\times\operatorname{SL}_2(q)}_{  m
\text{ \rm times}})\,:\, \operatorname{Sym}_m\le \operatorname{Sp}_n(q).
$$

Suppose that the action of $N_G(\Delta)$ on $\Delta$ is denoted by
$$
\rho:N_G(\Delta)\rightarrow \operatorname{Sym}(\Delta)\simeq
\operatorname{Sym}_m.
$$
By [21, Theorem~2],  $N_G(\Delta)^\rho=\operatorname{Sym}(\Delta).$  By Lemma~13
it follows that a $\pi$-Hall subgroup $H^\rho$ is either maximal in
$\operatorname{Sym}(\Delta)$ or equal to $\operatorname{Sym}(\Delta)$. In
particular,
$$
H^\rho\operatorname{prn}\operatorname{Sym}(\Delta)
\quad \text{ and }\quad N_{\operatorname{Sym}(\Delta)}(H^\rho)=H^\rho.
$$
Since
$N_G(S)\le N_G(\Delta)$ and $g\in N_G(S)$ there exists an element
$y\in\langle H,H^g\rangle$ such that $(H^g)^\rho=(H^y)^\rho$. So
$$
(gy^{-1})^\rho\in N_{\operatorname{Sym}(\Delta)}(H^\rho)=H^\rho.
$$

Denote by $A$ the kernel of $\rho$. The structure of $N_G(\Delta)$  implies
that if
$\overline{\phantom{X}}:A\rightarrow A/{O}_\infty(A)$ is a natural
homomorphism, then $\overline{A}$ possesses a normal subgroup isomorphic to
 $$
 \underbrace{\operatorname{PSL}_2(q)\times\dots
\times\operatorname{PSL}_2(q)}_{ m \text{ \rm times}},
 $$
and index of the subgroup in $\overline{A}$ is a $2$-power. By Lemmas~7--9
(we take the class of solvable groups as ${\frak X}$) and Lemma~25  we conclude
that $\pi$-Hall subgroups of $A$ are pronormal. Now by $\pi$-Hall subgroups of
$HA$ are pronormal by Lemma~8. Moreover, $gy^{-1}\in HA$ since
$(gy^{-1})^\rho\in H^\rho$. Therefore  $H^z=H^{gy^{-1}}$ for some
$z\in \langle H, H^{gy^{-1}}\rangle\le\langle H, H^{g}\rangle$.  Let $x=zy$.
Then  $H^x=H^g$ and $x\in \langle H, H^{g}\rangle$.
\end{proof}

Thus, if either statement~(6) of Lemma~12 holds, or statement~(7) of the same
lemma holds and for the preimage $H^*\le G^*=\operatorname{SL}^\eta(q)$ of
$H$ statement~(3) of Lemma~18 holds, then $H\operatorname{prn} G$. Notice also
that statements~(7) of Lemma~12 and ~(1) of Lemma~18 hold, then
$H\operatorname{prn} G$ by Lemma~25.

The next lemma allows to exclude also the case, when statements~(7) of Lemma~12
and (4) of Lemma~18 hold.

\begin{Lemma}
Let $G=\operatorname{PSL}_4^\eta(q)$, where $q$ is odd. Then $N_G(S)=S$.
\end{Lemma}

\begin{proof}
The claim follows by Lemma~12 since the $2$-adic expansion of $4$ has only one
unit.
\end{proof}

In case, when statements~(7) of Lemma~12 and~(2) of Lemma~18 hold,
$H$ normalizes a maximal torus of order
$(q-\eta)^{n-1}/(n,q-\eta)$ of $G=\operatorname{PSL}^\eta_n(q)$. We
consider this case as statement (5) of Lemma~12 in the next lemma.

\begin{Lemma}
Assume that $(\star)$ holds and
$G=\operatorname{PSL}_n^\eta(q)$, where $q$ is a power of a prime $p\not\in\pi$.
Suppose als that $q\equiv \eta \pmod 4$ and there exists a maximal
$H$-invariant torus $T$  of order
$(q-\eta)^{n-1}/(n,q-\eta)$.
Then

$(1)$~$N_G(T)=N(G,T)$;

$(2)$~$N_G(T)/T\simeq \operatorname{Sym}_n$;

$(3)$~$N_G(T)$ includes $N_G(S)$;

$(4)$~$H$ and $H^g$ are conjugate in $\langle H,H^g\rangle$.
\end{Lemma}

\begin{proof}
Statement (1) follows from [10, Lemma~3.10], since $T$ is invariant under given
Sylow  $2$-subgroup $S$ of $G$.

(2)~Since the identity $N_G(T)=N(G,T)$ holds, the factor group
$N_G(T)/T=N(G,T)/T$ is isomorphic to $\operatorname{Sym}_n$ (this factor group
is included in the Weyl group of $G$, which is isomorphic
to $\operatorname{Sym}_n$, on the other hand, a subgroup of type
$T.\operatorname{Sym}_n$ is included in $G$ ans so in $N_G(T)$).

(3)~Since $H\leq N_G(T)$ and $3\in\pi$, $N_G(T)$ includes a Sylow
$3$-subgroup of $G$. So, by [10, Lemma~3.13], it follows that
$N_G(S)\leq N_G(T)$.

(4)~In view of statement (3) of the lemma we remain to prove that
$H \operatorname{prn} N_G(T)$.  By statement~(2) of the lemma,
$N_G(T)$ is an extension of an abelian group $T$ by
$\operatorname{Sym}_n$. Consider the natural epimorphism
$\overline{\phantom{G}}:N_G(T) \rightarrow N_G(T)/T\simeq
\operatorname{Sym}_n$. By [20, Lemma~2.1(a)], $\overline{H}$ is a
$\pi$-Hall subgroup of $\operatorname{Sym}_n$.  Since, in view of the condition
 $2,3\in\pi$ and Lemma~13, each $\pi$-Hall subgroup is either maximal in
$\operatorname{Sym}_n$, or equal to $\operatorname{Sym}_n$, we have
$\overline{H}\operatorname{prn}\overline{N_G(T)}$. Taking the class of all
abelian groups as ${\frak X}$ in Lemma~9 we obtain that  ${H \operatorname{prn}
N_G(T)}$.
\end{proof}

Thus we have considered all possible cases, except the case, when statement
(7) of Lemma~12 holds,
$G=\operatorname{PSL}_{11}^\eta(q)$, and for the preimage
$H^*\le\operatorname{SL}_{11}^\eta(q)$ of $H$ statement (5) of Lemma~18 holds.
In particular the following lemma is true.

\begin{Lemma}
Let $2,3\in\pi$ and $q$ be a power of a prime $p\not\in\pi$. Then $\pi$-Hall
subgroups in $\operatorname{PSL}_n^\eta(q)$,
$\operatorname{PGL}_n^\eta(q)$,
$\operatorname{SL}_n^\eta(q)$, and $\operatorname{GL}_n^\eta(q)$ for
$n\leq 4$ and $n=8$ are pronormal.
\end{Lemma}

\begin{proof}
For $\operatorname{PSL}_n^\eta(q)$ the lemma follows directly from Lemma~18.
For $\operatorname{PGL}_n^\eta(q)$ the claim follows from Lemma~8 since
$$
\big|\operatorname{PGL}_n^\eta(q):\operatorname{PSL}_n^\eta(q)\big|=(n,q-\eta)
$$
divides $n$ and so it is a $\pi$-number. Finally, $\operatorname{SL}_n^\eta(q)$
and $\operatorname{GL}_n^\eta(q)$  are extensions of abelian groups by
$\operatorname{PSL}_n^\eta(q)$ and $\operatorname{PGL}_n^\eta(q)$. The
assertion of the lemma follows from above arguments and Lemma~9.
\end{proof}

Consider the remaining case. We need

\begin{Lemma}
Let $G^*=\operatorname{SL}_{11}^\eta(q)$, $q$ be odd, and
$S^*\in\operatorname{Syl}_2(G^*)$. Set $\Delta=\operatorname{Fun}_G(S^*)$.
Then

$(1)$~$|\Delta|=5$ and $S^*$ acting on $\Delta$ has exactly two orbits: $\Gamma$
of order~$4$ and  $\Gamma_0$ of order~$1$;

$(2)$~$\Gamma$ are $\Gamma_0$ $N_{G^*}(S^*)$-invariant;

$(3)$~if  $\Gamma'$ is an $S^*$-invariant set of pairwise commuting fundamental
subgroups of
$G^*$ such that $|\Gamma'|=4$, then $\Gamma'=\Gamma$.
\end{Lemma}

\begin{proof}
Denote by $\rho$ the action of $N_{G^*}(\Delta)$ on $\Delta$. According to [21,
Theorem 2]
$$
N_{G^*}(\Delta)^\rho=\operatorname{Sym}(\Delta)\simeq \operatorname{Sym}_5.
$$
$S^\rho$ is a Sylow $2$-subgroup of $\operatorname{Sym}_5$ and so it has two
orbits on $\Delta$: one orbit of length $4$ and another of length $1$.
This implies statement $(1)$. Statement $(2)$ follows from the fact that
$S^*$ and $\Delta$ are $N_{G^*}(S^*)$-invariant. Finally, $\Gamma'$ is included
in $\Delta$, since $\Delta$ is a unique maximal $S^*$-invariant set of pairwise
commuting fundamental subgroup. So $\Gamma'$ is a union of some orbits of
$S^*$ on $\Delta$ and, in view of $(1)$, equals~$\Gamma$.
\end{proof}

\begin{Lemma}
Let $G^*=\operatorname{SL}^\eta_{11}(q)$ be a special linear o unitary group
and  $V$ be its natural module equipped with a trivial or unitary form
respectively. Assume that $H^*\in\operatorname{Hall}_{\pi}(G^*)$, where
$\pi\cap\pi(G^*)=\{2,3\}$, and suppose that $H^*$ is included in a subgroup of
type
$$
L=\big(\big(\operatorname{GL}_2^\eta(q)\wr \operatorname{Sym}_4\big)
\times\big(\operatorname{GL}_1^\eta(q)\wr\operatorname{Sym}_3\big)\big)\cap G^*.
$$
Let
$S^*\in\operatorname{Syl}_2(H^*)\subseteq\operatorname{Syl}_2(G^*)$
and $g^*\in N_{G^*}(S^*)$. Then

$(1)$~$H^*$ leaves invariant a set $\Gamma'=\{K_1,K_2,K_3,K_4\}$ consisting
from pairwise commuting fundamental subgroups;

$(2)$~$\Gamma'$ is $N_{G^*}(S^*)$-invariant;

$(3)$~if $V_i=[K_i, V]$ and $U=\sum V_i$, then $U$ is invariant under both
$H^*$ and $N_{G^*}(S^*)$;

$(4)$~the stabilizer $M$ in $G^*$ of $U$ is a subgroup with pronormal
$\pi$-Hall subgroups;

$(5)$~$H^*\operatorname{prn} G^*$.
\end{Lemma}

\begin{proof}
Consider a subgroup $\big(\operatorname{GL}_2^\eta(q)\wr
\operatorname{Sym}_4\big)\cap G^*$ of
$L=\big(\big(\operatorname{GL}_2^\eta(q)\wr \operatorname{Sym}_4\big)
\times\big(\operatorname{GL}_1^\eta(q)\wr\operatorname{Sym}_3\big)\big)\cap
G^*$, and in the base of the wreath product consider distinct normal subgroups
$K_1, K_2,K_3,K_4$ isomorphic to $\operatorname{SL}_2(q)$. Clearly,
$K_i\not\in G^*$ and $K_i\not\in L$ for all $i=1,2,3,4$. Moreover, the set
$\Gamma'=\{K_1, K_2,K_3,K_4\}$ is $L$-invariant and so is  $H^*$-invariant.
Statement  $(1)$ is proven. Statement $(2)$ follows from Lemma~31.  Notice that
$V_i$ can be considered as the natural module for $K_i$, therefore $\dim
(V_i)=2$ and $V_i\cap V_j=0$ for $i\ne j$.  In particular, $\dim(U)=8$. Since
$\Gamma'$ is invariant under both $H^*$ and $N_{G^*}(S^*)$ it follows that
the set $\{ V_1,V_2,V_3,V_4\}$, and so the subspace $U$, are also invariant
under both $H^*$ and $N_{G^*}(S^*)$. Thus  $(3)$ is proven. If $\eta=+$, then
the stabilizer~$M$ of $U$ is an extension of a $p$-group by a central product
$\operatorname{GL}_{8}(q)\circ \operatorname{GL}_3(3)$
(see [11, Proposition 4.1.17]), and by Lemmas~30,~7, and~9  we conclude that
$\pi$-Hall subgroups of $M$ are pronormal. If $\eta=-$, then  $M$ is isomorphic
to a central product $\operatorname{GU}_8(q)\circ \operatorname{GU}_3(q)$
(see [11, Proposition 4.1.4]), and applying again Lemmas~30 and~7, we obtain
statement $(4)$. In view of $(3)$, $H^*$ and every $g^*\in N_{G^*}(S^*)$ are
included in $M$. Now from $(4)$ and Lemma~5 we conclude
that $H^*\operatorname{prn} G^*$.
\end{proof}

We continue the proof of the theorem and consider the remaining case. Assume
that statement~(7) of Lemma~12 holds,
$G=\operatorname{PSL}_{11}^\eta(q)$, and for the preimage
$H^*\le\operatorname{SL}_{11}^\eta(q)$ of $H$ statement~(5) of Lemma~18 holds.
By Lemma~32 we have
$H^*\operatorname{prn} \operatorname{SL}_{11}^\eta(q)$. Applying Lemma~6 we
conclude the proof of Theorem~1.
\qed

\section{Conclusion}

In connection with the proof of Theorem~1 we make a small note. The proof is
naturally divided into two cases. The first case, when a Hall subgroup $H$ of a
finite simple group $G$ has  odd order (equivalently, even index). The proof in
this case is reduced to application of Hall theorem [7, Theorem A1] (Lemma~10)
and Gross theorem [9, Theorem B] (Lemma~11). In the second case, when a Hall
subgroup $H$ has even order (equivalently, odd index), the technique is
absolutely different. We use the fact that $H$ includes a Sylow
$2$-subgroup $S$ of $G$, and so, by  Lemma~5, we need to check that
$H$ and $H^g$ are conjugate in $\langle H,H^g\rangle$ only for those
$g$,  that normalize $S$. Then we apply the structure of normalizers of Sylow
$2$-subgroups in finite simple groups obtained by A.~S.~Kondrat'ev
(Lemma~12). This technique could be probably applied in a more general
situation, For example, the following conjecture is of interest.

\begin{Conj}
Subgroups of odd index are pronormal in finite simple groups.
\end{Conj}

In view of Lemma~5, Conjecture~1 holds for all finite simple groups possessing
a self-normalizing Sylow  $2$-subgroup (for example, according to
Kondrat'ev theorem (Lemma~12) in alternating groups of degree greater than
$5$, in orthogonal groups, and in most classes of sporadic and exceptional
groups).

\end{document}